\documentclass[11pt]{article} 
\usepackage{amsfonts,amsmath,amssymb,color,amscd}
\usepackage[all,cmtip]{xy}
\usepackage[T1]{fontenc} 
\usepackage[francais,english]{babel} 
\usepackage{typearea}

\newcommand{\cqfd}{\hfill $\square$}
\newcommand\N{{\mathbb N}}

\newcommand\R{{\mathbb R}}
\newcommand\A{{\mathbb A}}
\newcommand\Z{{\mathbb Z}}
\newcommand\C{{\mathbb C}}

\newcommand\tr{\hbox to 1mm  {${}^t \!  $} }

\begin{document}

\centerline{\bf \Large A Characterization of Semisimple Plane Polynomial}

\centerline{\bf \Large Automorphisms.}

\vskip6mm

\centerline{ Jean-Philippe FURTER,}

\centerline{Dpt. of Math., Univ. of La Rochelle,}
 
\centerline{av. M. Cr\'epeau, 17 000 La Rochelle, FRANCE}

\centerline{email: jpfurter@univ-lr.fr}

\vskip5mm

\centerline{Stefan MAUBACH.\footnote{Funded by Veni-grant of council for the
physical sciences, Netherlands Organisation for scientific research (NWO). Partially funded by the Mathematisches
Forschungsinstitut Oberwolfach as an Oberwolfach-Leibniz-Fellow.}}

\centerline{ Dpt. of Math., Radboud University Nijmegen (RUN), }

\centerline{ postbus 9010, 6500 GL Nijmegen, THE NETHERLANDS }

\centerline{ email: s.maubach@science.ru.nl}

\vskip5mm

\centerline{April 10, 2008.}

\vskip10mm

\centerline{\bf Abstract.}
It is well-known that an element of the linear group
${\rm GL}_n( \C)$ is semisimple
if and only if its conjugacy class is Zariski closed.
The aim of this paper is to show that the same result holds for the
group of complex plane polynomial automorphisms.

\vskip4mm

\centerline{\bf Keywords.}

\vskip1mm

\centerline{Affine space, Polynomial automorphisms.}

\vskip10mm

\centerline{\bf  I. INTRODUCTION.}

\vskip5mm

If $K$ is any commutative ring, a polynomial endomorphism of the affine plane $\A^2_K$ over $K$
will be identified with its sequence
$f=(f_1,f_2)$ of coordinate functions $f_j \in K[X,Y]$.
We define the
degree of $f$ by
${\displaystyle \mbox{deg } f= \max_{1 \, \leq \,  j \,  \leq \, 2}\mbox{deg }f_j}$.

Let ${\cal G}$ be the group of polynomial automorphisms
of $\A^2_{\C}$
and let ${\cal G} (K)$ be  the group of polynomial automorphisms
of $\A^2_K$.

\vskip5mm

In linear algebra it is a well-known result that an element of ${\rm GL}_n(\C)$
has a closed conjugacy class if and only if it is semisimple, i.e. diagonalizable.
This is a very useful characterization, especially from a group action viewpoint.
It is a natural question to ask if a polynomial automorphism is semisimple
if and only if its conjugacy class is closed in the set of polynomial automorphisms.
This last statement hides two definitions: what is a semisimple polynomial automorphism
and what topology does one have on the group of polynomial automorphisms?

According to \cite{Fu-Ma},  the usual notion of  semisimplicity can be extended
from the linear to the polynomial case
by saying that a polynomial automorphism
is semisimple if it admits a vanishing polynomial with single roots.
In this paper we restrict to the dimension 2.
In this case, we will show below (see subsection 2.5) that it is equivalent to saying that
it is diagonalizable, i.e. conjugate to some diagonal automorphism
$(aX,bY)$ where $a,b \in \C^*$.

The topology of the group of polynomial automorphisms 
has been defined in \cite{Sha1,Sha2}. Let us describe it in
dimension 2 (the description would be analogous in dimension $n$).
The space ${\cal E}:=\C[X,Y]^2$  of polynomial endomorphisms
of $\A^2_{\C}$ is naturally an infinite dimensional algebraic variety
(see \cite{Sha1,Sha2} for the definition).
This roughly means that
${\cal E}_{\leq \, m}:= \{ f \in {\cal E}, \, {\rm deg}\, f \leq m \}$
is a (finite dimensional) algebraic variety
for any $m \geq 1$, which comes out from the fact that
it is  an affine space.
If $Z \subseteq {\cal E}$, we  set $Z_{\leq \, m} := Z \cap {\cal E}_{\leq \, m}$.
The space ${\cal E}$ is endowed with the
topology of the inductive limit, in which
$Z$ is closed (resp. open, resp. locally closed) if and only if
$Z_{\leq \, m}$ is
closed  (resp. open, resp. locally closed) in ${\cal E}_{\leq \, m}$ for any $m$.
Since ${\cal G}$ is locally closed in ${\cal E}$
(see \cite{B-C-W,Sha1,Sha2}),
it is naturally an infinite dimensional algebraic variety.

The aim of this paper is to show the following result.

\vskip5mm

\noindent {\bf Main Theorem.} A complex plane polynomial automorphism is
semisimple if and only if its conjugacy class is closed.

\vskip5mm

\noindent {\bf Application.}
If $f$ is a finite-order automorphism of the affine space $\A^3_{\C}$,
it is still unknown whether or not  it is diagonalizable.
Since any
finite-order linear automorphism is diagonalizable,
it amounts to saying that $f$ is linearizable, i.e. conjugate
to some linear automorphism.
To our knowledge,
even the case where $f$ fixes the last coordinate
is unsolved.
In this latter case, $f$ is traditionally seen as an element of 
${\cal G}( \C[Z] )$.
For each $z \in \C$, let $f_z \in {\cal G}$ be the automorphism induced by $f$
on the plane $Z=z$.
Using the amalgamated structure of  ${\cal G}( \C(Z) )$,
we know that $f$ is conjugate
in this group to some  $(aX,bY)$,
where $a,b \in  \C^*$ (see \cite{Ju,Ku,Se}).
This implies that $f_z$ is generically conjugate to $(aX,bY)$,
i.e. for all values of $z$ except perhaps finitely many.
The above theorem shows us that there is no exception:
for all $z$, $f_z$ is conjugate to $(aX,bY)$.
This could be one step for showing that such an $f$ is diagonalizable
in the group of polynomial automorphisms of $\A^3_{\C}$.
One can even wonder if the following is true.

\vskip5mm

\noindent {\bf Question 1.1.} Is  any finite-order automorphism
belonging to ${\cal G}( \C[Z] )$  diagonalizable in this group?

\vskip10mm

We begin in section 2  by studying
the so called locally finite plane polynomial automorphisms,
i.e. the automorphisms admitting a
non-zero vanishing polynomial.
The principal tool is the notion of pseudo-eigenvalues (see 2.2).
It is used for defining a trace (see 2.3)
and the subset ${\cal S} \subseteq {\cal G}$
of automorphisms admitting a single fixed
point (see 2.4).
Let us note that our text contains three natural questions 
which we were not yet able to answer.
Finally, we study the semisimple automorphisms
and show that their conjugacy class is characterized
by the pseudo-eigenvalues (see 2.5). 

The proof of the main theorem is given in section 3.
Subsection 3.1 is devoted to
an algebraic lemma whose proof relies on a valuative criterion
while subsection 3.2 is devoted to a few topological lemmas
(lemma 3.4 for example relies on Brouwer fixed point theorem).

\vskip10mm

\centerline{\bf  II. LOCALLY FINITE PLANE POLYNOMIAL AUTOMORPHISMS.}

\vskip5mm

\indent {\bf 1. Characterization.}

\vskip5mm

According to \cite{Fu-Ma}, a polynomial endomorphism is called locally finite
(LF for short) if it admits a non-zero vanishing polynomial.
The class of LF plane polynomial automorphisms will be denoted by ${\cal LF}$.
We recall that an automorphism is said to be triangularizable if it is
conjugate to some triangular automorphism $(aX+p(Y), bY+c)$, where $a,b \in \C^*$,
$c \in \C$ and $p \in \C[Y]$.
Using the amalgamated structure of ${\cal G}$, one can  show
the following:

\vskip5mm

\noindent {\bf Theorem 2.1.} If $f \in {\cal G}$, the following assertions are equivalent:

(i) $f$ is triangularizable;

(ii) the dynamical degree
$dd(f):= {\displaystyle \lim_{n \to \infty} } ({\rm deg}\,f^n)^{1/n}$
is equal to $1$;

\vspace{-1mm}

(iii) $\mbox{deg }f^2 ~ \leq ~ \mbox{deg }f$;

(iv) $\forall~ n \in \N,~  \mbox{deg }f^n ~ \leq ~ \mbox{deg }f$;

(v)  for each   $\xi \in \A^2$, the sequence $n \mapsto f^n( \xi)$ is a linear recurrence sequence;

(vi) $f$ is LF.

\vskip5mm

\noindent {\bf Proof.} For (i-ii), (iii-iv),  (v) and (vi) ,
see respectively \cite{F-M}, \cite{Fu1}, \cite{Fu2} and \cite{Fu-Ma}. \cqfd

\vskip5mm

In this case, the minimal polynomial $\mu_f$ of $f$ is defined as the (unique)
monic polynomial generating the ideal
$\{ p \in \C[T], \, p(f) = 0 \}$.
Let us note that even if the class ${\cal LF}$ is invariant by conjugation,
the minimal polynomial is not.

\vskip5mm

\noindent {\bf Corollary 2.1.} ${\cal LF}$ is closed in ${\cal G}$.

\vskip5mm

\noindent {\bf Proof.}
By assertions (ii, iv),
we have   ${\cal LF}_{\leq \, m}=
\{ f \in {\cal G}, \, \forall \, n \in \N, \,
{\rm deg}\,f^n \, \leq m \}$ (for any $m \geq 1$).
This proves that 
${\cal LF}_{\leq \, m}$ is closed in ${\cal G}_{\leq \, m}$. \cqfd

\vskip5mm

\indent {\bf 2. The pseudo-eigenvalues.}

\vskip5mm

If $f \in {\cal LF}$, it is conjugate to some
triangular automorphism $t=(aX+p(Y), bY+c)$.
It is explained in \cite{F-M} (cf. the remark  on page 87) 
that the unordered pair $\{ a,b \}$ is an invariant:
if $t$ has a fixed point, then $a$ and $b$ are equal to
the two eigenvalues of the derivative at that fixed point
and if $t$ has no fixed point, then  the pair $\{ a,b \}$
must be equal to $\{ 1, {\rm Jac}\,f \}$.

\vskip5mm

\noindent {\bf Definition.} $a,b$ are called the pseudo-eigenvalues of $f$.

\vskip5mm

Let $<a,b> \, := \{a^kb^l, \, k,l \in \N \}$
be the submonoid of $\C^*$ generated by $a,b$ and
let $f^* : \, r \mapsto r \circ f$
be the algebra automorphism of $\C[X,Y]$ associated to $f$.
The following result relates the
pseudo-eigenvalues of $f$
with the eigenvalues of $f^*$.

\vskip5mm

\noindent {\bf Lemma 2.1.} If $a, b $ are the pseudo-eigenvalues
of  $f \in {\cal LF}$, then
$<a,b>$ is the set of eigenvalues of $f^*$.

\vskip5mm

\noindent {\bf Proof.}
We may assume that $f=(aX+p(Y), bY+c)$.
Let $d$ be the degree of $p(Y)$.

Let $M:=\{ X^kY^l, \, k,l \geq 0 \}$ be the set of
all monomials in $X,Y$ and let us endow $M$ with
the monomial order $\prec$ (see  \cite{E}) defined by

\centerline{$X^kY^l \prec X^m Y^n$ $\Longleftrightarrow$ $k<m$ or ($k=m$ and $l<n$). }

For any $s \geq 0$, we observe
that the vector space $V_s$ generated by
the $X^kY^l$ such that $dk+l \leq s$ is stable by $f^*$.
Let us denote by $f^*_{||V_s}$ the induced linear endomorphism of $V_s$.

Since $f^*(X^kY^l)-a^kb^l X^kY^l \in {\rm Span}(X^mY^n)_{X^mY^n \, \prec \, X^kY^l}$
(exercise), the matrix of $f^*_{||V_s}$ in the basis
$X^kY^l$ (where the $X^kY^l$ are taken with the order $\prec$) is upper
triangular with the $a^kb^l$'s on the diagonal.
The result follows from the equality
$\C[X,Y]= {\displaystyle \bigcup_s } \, V_s$. \cqfd

\vskip5mm

It is well-known that the eigenvalues of a linear automorphism are roots of its minimal polynomial.
The same result holds for LF plane polynomial automorphisms:

\vskip5mm

\noindent {\bf Lemma 2.2.} The pseudo-eigenvalues  are roots of the minimal polynomial.

\vskip5mm

\noindent {\bf Proof.} We will use the basic language of
linear recurrence sequence that we now recall
(see  \cite{CMP} for details).
If $U$ is any complex vector space, the set of
sequences $u~:~ \N \to U$ will be denoted by $U^{\N}$. For
$p= {\displaystyle \sum_kp_k~T^k \in \C[T]}$, we define $p(u)
\in U^{\N}$ by the formula

\vskip2mm

\centerline{${\displaystyle \forall~n \in \N,~\left( \rule{0mm}{3.5mm} p(u)
\right) (n)= \sum_kp_k~u(n+k).}$}

\vskip2mm

Let $U[T]$ be  the set of polynomials in $T$
with coefficients in $U$, alias the set of polynomial maps from
$\C$ to $U$.

\vskip3mm

The theory of linear recurrence sequence
relies on the fact that if 
${\displaystyle p=\alpha ~\prod_{1 \, \leq \, k \, \leq \, c} \! (T-
\omega_k)^{r_k} }$ is the decomposition into irreducible factors
of some non-zero polynomial $p$, then
$p(u)=0$ if and only if there exist $q_1,\ldots,q_c \in U[T]$ with
$\mbox{deg }q_k~ \leq r_k -1$ such that

\vskip2mm

\centerline{${\displaystyle \forall~n \in \N ,~ u(n)~=~\sum_{1 \leq k \leq c} \! \omega_k^nq_k(n) }$.}

\vskip2mm

As a consequence, it is clear that
${\cal I}_u:=\{ p \in \C[T],~ p(u)=0 \}$ is an ideal of $\C[T]$.
We say that $u$ is  a linear recurrence sequence when
${\cal I}_u \neq \{ 0 \}$.
In this case, the minimal polynomial of $u$ is the (unique)
monic polynomial $\mu _u$ generating the ideal ${\cal I}_u$.

We say that $u$ is of exponential type if the following equivalent  assertions
are satisfied:

\vskip2mm

(i) there exist $\omega _1, \ldots, \omega_c \in \C$, $q_1,\ldots,q_c \in U$ such that
${\displaystyle \forall~n ,~ u(n)~=~\sum_{1 \leq k \leq c} \! \omega_k^nq_k }$.

\vskip2mm

(ii) $\mu_f$ has single roots.

\vskip2mm

If $l:\, U \to V$ is any linear map,
let us note that $v:=l(u)$ is still a linear recurrence
sequence and that $\mu_v$ divides $\mu_u$.

\vskip2mm

If $A \in M_k(\C)$ is a square matrix, one could easily check
that the minimal polynomial of $A$ is equal to the minimal
polynomial of the linear recurrence sequence $n \mapsto A^n$.

\vskip1mm

Let now  $f \in {\cal LF}$ be a LF plane polynomial automorphism.
One could also check
that the minimal polynomial of $f$
is equal to the minimal polynomial of the linear recurrence
sequence $n \mapsto f^n$ (see \cite{Fu-Ma} for details).

\vskip2mm

Let us now begin the proof.

\vskip4mm

\noindent \underline{First case.} $f$ admits at least one fixed point $\xi$.

\noindent If  $(\omega_i)_{1 \, \leq \,  i \,  \leq \, r}$ are the roots of $\mu_f$, there exist polynomial
endomorphisms $h_{i,j}$ such that
$f^n= {\displaystyle \sum_{i,j} } \omega_i^n n ^j h_{i,j}$ for any $n \geq 0$.
Differentiating at the point $\xi$, we get
$f'(\xi)^n= {\displaystyle \sum_{i,j} } \omega_i^n n ^j (h_{i,j})'(\xi)$,
so that the eigenvalues of the matrix $f'(\xi)$
are among the $\omega_i$'s.
But since $\xi$ is a fixed point, 
these eigenvalues are
the pseudo-eigenvalues of $f$. 

\vskip4mm

\noindent \underline{Second case.} $f$ admits no fixed point.

\noindent By  theorem 3.5 of \cite{F-M}, $f$ can be expressed as
$f=\varphi \circ t \circ \varphi^{-1}$
where $\varphi \in {\cal G}$, $p \in \C [ Y]$, $b \in \C^*$ and either

(i) $t= (X+1,bY)$;

(ii) $t= (X+p(Y^r),bY)$
where $r \geq 2$, $b^r=1$, $p(0)=1$;

(iii) $t= (X+p(Y),Y)$.

\vskip4mm

\indent \underline{Subcase (i).}

We have $f^n= \varphi \circ (X+n,b^nY) \circ \varphi ^{-1}$ for any $n \geq 0$.
Let $\psi := \varphi^{-1}$ and let $(e_1,e_2)$ be the canonical basis
of the $\C[X,Y]$-module $\C[X,Y]^2$.
Since the family $ \psi_1^i \psi_2^j$ for $i,j \geq 0$ is a basis of  $\C[X,Y]$,
the family $ \psi_1^i \psi_2^j e_k$ is a basis of ${\cal E}= \C[X,Y]^2$.

If $\varphi_k= {\displaystyle \sum_{i,j} }\varphi_{k,i,j} \, X^iY^j$ for $k=1,2$, an easy computation
would show that the $\psi_1 e_k$-component of $f^n$ is
${\displaystyle \sum_i} i \varphi_{k, i,0} \, n^{i-1}$
and that the $\psi_2 e_k$-component of $f^n$ is 
${\displaystyle \sum_i }  \varphi_{k, i,1} \, n^{i}b^n$.
But the matrix
$\left[
\begin{array}{cc}
\varphi_{1,1,0} & \varphi_{1,0,1} \\
\varphi_{2,1,0} & \varphi_{2,0,1}
\end{array}
\right]$  corresponds to the linear part
of $\varphi$ so that it is invertible.
Therefore at least one of the $\varphi_{k,1,0}$ is non-zero
showing that $1$ is a root of the minimal polynomial
of the linear recurrence sequence sending $n$ to
the $\psi_1 e_k$-component of $f^n$.
Consequently, $1$ is a root 
of the linear recurrence sequence sending $n$ to $f^n$.
This means  that $\mu_f(1)=0$.
In the same way, at least one of the $\varphi_{k,0,1}$ is non-zero
showing that $\mu_f(b)=0$.

\vskip4mm

\indent \underline{Subcase (ii).}

We have $f^n= \varphi \circ (X+np(Y^r),b^nY) \circ \varphi ^{-1}$ for any $n \geq 0$.
We go on  as in subcase (i). The computations are slightly different,
but the results (and conclusions) are exactly the same.

\vskip4mm

\indent \underline{Subcase (iii).}

We recall that a linear recurrence sequence is polynomial
if and only if its minimal polynomial is of the kind $(T-1)^d$.
We conclude by noting that the sequence
$n \mapsto f^n$ is obviously polynomial. \cqfd

\vskip5mm

\indent {\bf 3. The trace.}

\vskip5mm

It is natural to set the following

\vskip5mm

\noindent {\bf Definition.} If $f \in {\cal LF}$ has pseudo-eigenvalues $\{ a,b \}$,
its trace is ${\rm Tr}\,f:=a+b$.

\vskip5mm

\noindent {\bf Remark.} The trace is by construction an invariant of
conjugation.
It is well-known that the Jacobian map
${\rm Jac}:\, {\cal G} \to \C^*$ 
also.
In the locally finite case, we have of course
${\rm Jac}\,f= ab$.

\vskip5mm

\noindent {\bf Question 2.1.} Is the map ${\rm Tr}:\, {\cal LF} \to \C$ regular?

\vskip5mm

This means that for any $m$ the restricted map
${\cal LF}_{\leq \, m} \to \C$ is regular.
This regularity  would imply
a positive answer  to the following

\vskip5mm

\noindent {\bf Question 2.2.} Is the map ${\rm Tr}:\, {\cal LF}_{\leq \, m} \to \C$
continuous for the transcendental topology?

\vskip5mm

\noindent {\bf Remark.} This continuity would easily prove
the most difficult point of our main theorem.
If $f,g$ are semisimple automorphisms such
that $g$ belongs to the closure of the conjugacy class of $f$,
we want to show that they have the same pseudo-eigenvalues.
Indeed, it is clear that ${\rm Jac}\,f={\rm Jac}\,g$ and 
the above continuity would show that ${\rm Tr}\,f={\rm Tr}\,g$.

\vskip5mm

\noindent {\bf Definition.} Let  ${\cal U}$ (resp. ${\cal S}$) be
the set of LF polynomial automorphisms whose pseudo-eigenvalues
are equal to $1$ (resp. are different from $1$).

\vskip5mm

\noindent {\bf Remarks.} 1.  By theorem 2.3 of \cite{Fu-Ma}
${\cal U}$ is the set of polynomial automorphisms $f$
satisfying the following equivalent assertions:

(i) $f$ is  unipotent,
i.e. $f$ is  annihilated by  $(T-1)^d$ for some $d$;

(ii) $f$ is  the exponential
of some locally nilpotent derivation of $\C [X,Y]$.

\vskip2mm

2. It is easy to check  that
${\cal S}$ is the set of LF automorphisms admitting a single  fixed point
(in fact, we will see in proposition 2.1 below that we can get rid of the LF hypothesis).

\vskip2mm

3. Since ${\cal U}= {\rm Tr}^{-1}(\{ 2 \}) \cap {\rm Jac}^{-1}( \{ 1 \})$
and ${\cal S}=\{ f \in {\cal LF}, \, {\rm Tr}(f) \neq 1 + {\rm Jac}(f) \}$,
the regularity of the trace  would imply directly that
${\cal U}$ (resp. ${\cal S}$) is closed (resp. open) in ${\cal LF}$.

\vskip4mm

Let us check that ${\cal U}$ is closed.
If $m \geq 1$, let $d$ be the dimension of 
${\cal E}_{\leq \, m}$ and let $p(T) =(T-1)^d \in \C [T]$.
By assertion (iv) of theorem 2.1, we get
${\cal U}_{\leq \, m}=\{ f \in {\cal E}_{\leq \, m}, \, p(f)=0 \}$.
This shows that ${\cal U}_{\leq \, m}$ is closed in ${\cal E}_{\leq \, m}$
for any $m$, i.e.  ${\cal U}$ is closed in ${\cal E}$.

\vskip5mm

We will show in the next subsection that ${\cal S}$ is open in ${\cal LF}$.

\vskip5mm

\indent {\bf 4. The set ${\cal S}$.}

\vskip5mm

\noindent {\bf Definition.} If $f,g$ are polynomial endomorphisms of $\A^2_{\C}$,
let us define their coincidence ideal  $\Delta(f,g)$ 
as the  ideal generated by the $f^*(p)-g^*(p)$, where
$p$ describes $\C[X,Y]$.

The coincidence ideal $\Delta (f , {\rm id})$
will be called the fixed point ideal of $f$.

\vskip5mm

\noindent {\bf Remarks.} 1. The closed points of ${\rm Spec}\, \C[X,Y]/\Delta (f,g)$
correspond to the points $\xi \in \A^2_{\C}$ such that $f(\xi)=g(\xi )$.

2. Using the relation
$f^*(uv)-g^*(uv)= f^*(u) [ f^*(v)-g^*(v)]+ g^*(v) [ f^*(u)- g^*(u)]$,
we see that if  the algebra $\C[X,Y]$ is generated by
the $u_k$ ($1 \leq k \leq l$),
then the ideal $\Delta (f,g)$ is generated by the
$f^*(u_k)-g^*(u_k)$ ($1 \leq k \leq l$).

3. In particular, $\Delta (f,g)=
\left( \rule{0mm}{4mm} f^*(X)-g^*(X), f^*(Y)-g^*(Y) \right)
=(f_1-g_1,f_2-g_2)$.

\vskip5mm

The computation of the set of fixed points of a
triangular automorphism is easy and left to the reader.
We obtain the following result
(see also lemma 3.8 of \cite{F-M}).

\vskip5mm

\noindent {\bf Lemma 2.3.} If $f \in {\cal LF}$,
the set of its fixed points   is either empty,
either a point of multiplicity $1$ (if $ f \in {\cal S}$)
or either a finite disjoint union of  subvarieties
isomorphic to $\A^1$.

\vskip5mm

Let us note that saying that an automorphism admits exactly 1 fixed point
with multiplicity $1$ amounts to saying
that its fixed point ideal is a maximal ideal of $\C[X,Y]$.
Using the amalgamated structure of ${\cal G}$, it is observed in \cite{F-M}
that a polynomial automorphism $f \in {\cal G}$ is either triangularizable
(i.e. belongs to ${\cal LF}$) or conjugate to
some cyclically reduced element $g$
(see I.1.3 in \cite{Se}  or page 70 in \cite{F-M} for the definition).
In this latter case, the degree $d$ of $g$ is $\geq 2$
and it is shown in theorem 3.1 of \cite{F-M} that
${\rm dim}\,  \C[X,Y] / \Delta (g, {\rm id}) \, = \, d$.
As a conclusion, we obtain the  nice characterization of elements of ${\cal S}$:

\vskip5mm

\noindent {\bf Proposition 2.1.} If $f \in {\cal G}$, the following assertions are
equivalent:

(i) $f \in {\cal S}$;

(ii) $f$ has a unique fixed point of multiplicity $1$;

(iii) the fixed point ideal of $f$ is a maximal ideal of $\C [X,Y]$.

\vskip5mm

The next result is taken from lemma 4.1 of \cite{Fu-Ma} 
and will be used to prove  propositions 2.2 and 2.3 below.

\vskip5mm

\noindent {\bf Lemma 2.4.} Any triangularizable automorphism  $f$ can  be triangularized
in a "good" way with respect to the degree: 
there exist a triangular automorphism
$t$ and an automorphism $\varphi$ such that
$f = \varphi \, \circ \, t \, \circ \,  \varphi^{-1}$ with
$\mbox{deg }f ~=~ \mbox{deg} \, t ~( \mbox{deg }\varphi)^2$. 

\vskip5mm

The vector space $\A^2_{\C}$ will be endowed
with the norm
$\parallel ( \alpha, \beta) \parallel = \sqrt{|\alpha |^2 + | \beta | ^2}$.
The open (resp. closed) ball of radius $R \geq 0$
centered at a point $\xi \in \A^2_{\C}$ will be denoted
by $B_{\xi,R}$ (resp. $B'_{\xi,R}$).
If $\xi = 0$, we will write $B_R$ (resp. $B'_R$) instead of $B_{0,R}$
(resp. $B'_{0,R}$).

Since ${\cal E}$ is composed of $C^{ \infty}$ maps
from $\A^2_{\C}$ to $\A^2_{\C}$,
it is  endowed  with the $C^k$-topology (for each $k \geq 0$)
which is the topology of uniform convergence of the $k$ first derivatives
on all compact subsets.
However,  ${\cal E}_{\leq \, m}$ being a finite dimensional complex vector space,
it admits a unique Hausdorff topological vector space structure.
Therefore, the $C^k$-topology on ${\cal E}_{\leq \, m}$ is nothing else than the
transcendental topology.
We finish these topological remarks by recalling that 
for any constructible subset of some complex algebraic variety,
the (Zariski-)closure  coincide with the transcendental closure
(see for example \cite{Mu}).

\vskip5mm

\noindent {\bf Proposition 2.2.} ${\cal S}$ is an open subset of ${\cal LF}$.

\vskip5mm

\noindent {\bf Proof.} We want to show that
${\cal S}_{\leq \, m}$ is open in  ${\cal LF}_{\leq \, m}$.

\vskip2mm

\underline{Claim.} ${\cal S}_{\leq \, m}$
is a constructible subset of ${\cal LF}_{\leq \, m}$.

Let ${\cal T}$ be  the variety of triangular automorphisms
of the form $(aX+p(Y), bY+c)$ where
$a,b \in \C \setminus \{ 0,1 \}$, $c \in \C$ and
$p \in \C[Y]$ is a polynomial of degree $\leq m$.

The image $W$  of the morphism
${\cal G}_{\leq \, m} \times {\cal T} \to {\cal G}$,
$(\varphi ,t) \mapsto \varphi \circ t \circ \varphi ^{-1}$
is constructible and ${\cal S}_{\leq \, m} = W \cap {\cal LF}_{\leq \, m}$
by lemma 2.4 so that the claim is proved.

\vskip2mm

It is enough to show that  ${\cal S}_{\leq \, m}$
is open for the transcendental topology.
Let $f$ be a given element of ${\cal S}_{\leq \, m}$
and let $\xi \in \A^2$ be its fixed point.
The map $F:= f- {\rm id}$ is a local diffeomorphism
near $\xi$ since  $F'(\xi)$ is invertible.
Let $\varepsilon, \eta > 0$ be such that
$B_{\eta} \subseteq F( B_{ \xi , \varepsilon} )$
and $\forall \, x \in B_{ \xi, \varepsilon}, | \det F'(x) | \geq \eta$.
If $g$ is "near" $f$ for the $C^1$-topology, then $G:=g-{\rm id}$
will be "near" $F$ so that we will have
$B_{\eta /2} \subseteq G( B_{ \xi, \varepsilon} )$
and $\forall \, x \in B_{ \xi, \varepsilon}, | \det G'(x) | \geq \eta /2$.
Therefore, $g$ will have an isolated fixed  point in $B_{ \xi, \varepsilon}$.
If  $g \in {\cal LF}$,
lemma 2.3 shows us
that $g \in {\cal S}$. \cqfd

\vskip5mm

The next statement is given on page 49 of \cite{He} (cf. the application of theorem 3).
The result is also given  for the  field of rationals
on page 312 of \cite{M-M}.
However, the proof   remains unchanged for the field of complex numbers.
Finally,  \S 57 of \cite{Sei} contains a similar result.

\vskip5mm

\noindent {\bf Theorem 2.2.} Let $K:=d + (sd)^{2^n}$.
If $p,p_1,\ldots,p_s \in \C[X_1,\ldots,X_n ]$ are of degree
$\leq d$ and if $p \in (p_1,\ldots, p_s)$,
there exist $\lambda_1, \ldots, \lambda_s \in \C[X_1,\ldots,X_n]$
such that

\vskip2mm

(i) $p = {\displaystyle \sum_{1 \, \leq \, i \, \leq \, s} } \lambda_i \, p_i$
\hspace{5mm} and \hspace{5mm} 
(ii) ${\rm deg} \, \lambda_i \leq K$ for all $i$.

\vskip5mm

If $f \in {\cal S}$, its fixed point $\xi=(\alpha, \beta ) \in \A^2$ is implicitely defined by
the equality of the ideals $(f_1-X,f_2-Y)$ and $(X- \alpha, Y- \beta )$.
Using theorem 2.2, one can  express more "effectively" $\alpha,\beta$
in terms of $f_1,f_2$. Indeed, if $m\geq 1$ and $K_m:=m+(2m)^4$,
then for any $f \in {\cal S}_{\leq \, m}$
there exist $\lambda_1, \ldots, \lambda_4 \in \C[X,Y]$
of degree $\leq K_m$ such that
$X- \alpha  = \lambda _1 \, (f_1-X) +  \lambda _2 \, (f_2-Y)$
and $Y- \beta  = \lambda _3 \,  (f_1-X) +  \lambda_4 \, (f_2-Y)$.
Even with such "effective" results, we were not able to
answer the following

\vskip5mm

\noindent {\bf Question 2.3.} Is  the map
$Fix: {\cal S} \to \A^2$ sending
$f \in {\cal S}$ to its unique fixed point regular?

\vskip5mm

This means that for any  $m$ the restricted map 
$ {\cal S}_{\leq \, m} \to \A^2$ is regular.
The proof of proposition 2.2 shows us at least that it is 
continuous for the transcendental topology.

\vskip5mm

\indent {\bf 5. Semisimple  automorphisms.}

\vskip5mm

According to \cite{Fu-Ma}, a plane polynomial automorphism $f$ is said to be semisimple 
if the following equivalent assertions are satisfied:

(i) $f^*$ is semisimple (i.e.  $\C[X,Y]$
admits a basis of eigenvectors);

(ii) $f \in {\cal LF}$ and $\mu_f$ has single roots;

(iii) $f$ admits a vanishing polynomial with single roots.

\vskip1mm

Let us note that the class of semisimple automorphisms
is invariant by conjugation.
Therefore,
it  results from proposition 2.3 below that
  (i-iii)  are still equivalent to:

\vskip1mm

(iv) $f$ is diagonalizable.

\vskip5mm

\noindent {\bf Lemma 2.5.} If $t=(aX+p(Y), bY+c)$ is a triangular
semisimple automorphism,  there exists
a triangular automorphism $\chi$ of the same degree  such that
$t = \chi \circ  (aX , bY) \circ \chi ^{-1}$.

\vskip5mm

\noindent {\bf Proof.}

\noindent \underline{First step.} Reduction to the case $c=0$.

If $b=1$, let us show  that $c=0$.
The second coordinate of the $n$-th iterate $t^n$ is $Y+nc$.
Since $t$ is semisimple, the sequence $n \mapsto Y+nc$ must be of exponential type
showing that $c=0$.

\vskip2mm

If $b \neq 1$, set $l:=(X,Y+\frac{c}{b-1})$ and  replace
$t$ by $l \circ t \circ l^{-1}=(aX+p(Y), bY)$.

\vskip2mm

\noindent \underline{Second step.} Reduction to the case $p=0$.

If $\chi := (X+q(Y),Y)$, we get
$\chi \circ  (aX , bY) \circ \chi ^{-1}=(aX+q(bY)-aq(Y),bY)$.
Let us write $p={ \sum _k} p_k Y^k$.
To show the existence of $q$ (of the same degree as $p$)
satisfying $q(bY)-aq(Y)=p(Y)$ it is enough to show that
$a=b^k$ implies $p_k=0$.

For any $n \geq 0$, let $u_n$ be the $Y^k$-coefficient of the first component of $t^n$.
If $a=b^k$, we get $u_{n+1}=au_n + p_ka^n$, so that
$u_n= na^{n-1}p_k$. The sequence $n \mapsto u_n$ being of exponential type,
we obtain $p_k=0$. \cqfd

\vskip5mm

Combining lemmas 2.4 and 2.5, any semisimple automorphism  can be written

\vskip1mm

\centerline{ $f = ( \varphi \circ \chi ) \circ (aX,bY) \circ ( \varphi \circ \chi )^{-1}$
with $\mbox{deg }f ~=~ \mbox{deg} \, \chi ~( \mbox{deg}\,\varphi)^2$.}

\vskip1mm

Since $\mbox{deg}( \varphi \circ \chi ) \leq \mbox{deg} \, \varphi ~ \mbox{deg} \, \chi
\leq \mbox{deg} \,f$, we get:

\vskip5mm

\noindent {\bf Proposition 2.3.} Any semisimple automorphism $f$ can be written
$f = \psi  \circ (aX,bY) \circ \psi ^{-1}$ where $\psi$ is an automorphism
satisfying $\mbox{deg} \,\psi  \, \leq \, \mbox{deg} \,f$.

\vskip5mm

\noindent {\bf Corollary 2.2.} Two semisimple automorphisms are conjugate
if and only if they have the same pseudo-eigenvalues.

\vskip5mm

If $f \in {\cal G}$, let ${\cal C}(f):= 
\{ \varphi \circ f  \circ \varphi ^{-1}, \, \varphi  \in {\cal G} \}$
be its conjugacy class.
By definition, ${\cal C}(f)$ is closed in ${\cal G}$ if and only if
${\cal C}(f)_{\leq \, m}$ is closed in ${\cal G}$ for any $m \geq 1$.
However, if $Z \subseteq {\cal G}$, let us note that in general,
we do not have
$\overline{Z} = {\displaystyle \bigcup_{m \,  \geq \,  1} \overline{Z_{\leq \, m} } }$.

\vskip5mm

\noindent {\bf Corollary 2.3.} If $f$ is a semisimple automorphism,
then ${\cal C}(f)_{\leq \, m}$
is a constructible subset of ${\cal E}_{\leq \, m}$ (for any $m \geq 1$).

\vskip5mm

\noindent {\bf Proof.} We can assume that $f=(aX,bY)$.
The image $Z$ of the map
${\cal G}_{\leq \, m} \to {\cal G}$,
$\varphi \mapsto  \varphi \circ f \circ \varphi ^{-1}$ is constructible
and  ${\cal C}(f)_{\leq \, m}= Z \cap {\cal G}_{\leq \, m}$ by proposition 2.3. \cqfd

\vskip5mm

\noindent {\bf Remarks.}  1. This result shows us that the Zariski-closure of
${\cal C}(f)_{\leq \, m}$ coincide with its transcendental closure
(see subsection 3.2).

2. One could show that ${\cal C}(f)_{\leq \, m}$ is a constructible
subset of ${\cal E}_{\leq \, m}$ for any $f$, but we do not need this result.

\vskip5mm

\noindent {\bf Lemma 2.6.} If $f$ is semisimple, 
any element of  $\overline{ {\cal C}(f)_{\leq \, m}}$  also.

\vskip5mm

\noindent {\bf Proof.} We may assume that $f=(aX,bY)$.
Any element  which is linearly conjugate to $f$ is annihilated by $\mu _f$,
but for a general element of ${\cal C}(f)$,
this  is no longer true.
However,  we will build a polynomial $p$ with single roots
 annihilating any element of ${\cal C}(f)_{\leq \, m}$.
By propostion 2.3, any $g \in {\cal C}(f)_{\leq \, m}$
can be written $g= \varphi \circ f \circ \varphi^{-1}$
with ${\rm deg}\, \varphi \leq m$.
Therefore, for any $n \geq 0$, we have
$g^n= \varphi \circ (a^nX,b^nY) \circ \varphi^{-1}$.
If we set $\Omega:= \{ a^k b^l, \, 0 \leq k+l \leq m \}$,
there exists a family of polynomial endomorphisms $h_{\omega}$ ($\omega \in \Omega$)
such that  $g^n = {\displaystyle \sum_{ \omega \in \Omega} } \omega ^n h_{\omega}$
for any $n$.
In other words (see \S 1.3 of \cite{Fu-Ma} for details),  $p(g)=0$, where
$p(T):= {\displaystyle \prod_{\omega \in \Omega} } (T - \omega)$.
The equality $p(g)=0$ remains true if
$g \in \overline{ {\cal C}(f)_{\leq \, m}}$. \cqfd

\vskip10mm

\centerline{\bf  III. PROOF OF THE MAIN THEOREM.}

\vskip5mm

{\bf 1. Algebraic lemma.}

\vskip5mm

The  aim of this subsection is to prove the following
result which in some sense means  that
the spectrum of a linear endomorphism
remains unchanged  at the limit (see lemma 2.1).

\vskip5mm

\noindent {\bf Lemma 3.1.} Let $f=(aX,bY) \in {\cal G}$.
If $( \alpha X, \beta Y) \in \overline{{\cal C}(f)_{\leq \, m} }$, then
$ < \alpha , \beta > \,  = \,  <a,b> $.

\vskip5mm

Our proof will use a valuative criterion that we  give below.
We are indebted to Michel Brion for his useful advice on this subject.
Even if such a criterion sounds familiar
(see for example \cite{M-F-K}, chap. 2, \S\hspace{0.5mm}1, pp 52-54
or \cite{Gr}, \S\hspace{0.5mm}7), we have given a brief proof of it
for the sake of completeness.

Let $\C [[t]]$ be the algebra of complex formal power series
and let $\C ((t))$ be its quotient field.
If $V$ is a complex algebraic variety and
$A$ a complex algebra, $V(A)$ will denote
the  points of $V$ with values in $A$,
i.e. the set of morphisms ${\rm Spec}\, A \to V$.
If  $v$ is a closed point of $V$ and $\varphi  \in V \! \!  \left( \rule{0mm}{4mm} \C((t)) \right)$,
we will write
${\displaystyle v = \lim_{t \to 0} \varphi (t) }$ when:

(i) the point $\varphi : {\rm Spec}\, \C((t)) \to V$ is a composition
${\rm Spec}\, \C((t)) \to {\rm Spec}\, \C[[t]] \to V$;

(ii) $v$ is the point ${\rm Spec}\, \C \to {\rm Spec}\, \C[[t]] \to V$.

\vskip3mm

For example, if
$V=\A^1_{\C}$ and $ \varphi  \in V \! \!  \left( \rule{0mm}{4mm} \C((t)) \right)= \C((t))$,
we will write ${\displaystyle v = \lim_{t \to 0} \varphi (t) }$ when $\varphi  \in \C[[t]]$
and $v= \varphi(0)$.

\vskip3mm

\vskip5mm

\noindent {\bf Valuative criterion.} Let $f:~V \to W$ be a morphism of
complex algebraic varieties and let $w$ be a closed point of $W$.
The two following assertions are equivalent:

(i) $w \in \overline{ f(V) }$;

(ii) ${\displaystyle w= \lim_{t \to 0} f( \varphi (t) ) }$ for some 
$\varphi  \in V \! \! \left( \rule{0mm}{4mm}  \C ((t)) \right)$.

\vskip5mm

\noindent {\bf Proof.}

(i) $\Longrightarrow$ (ii). If $w \in \overline{ f(V) } \setminus f(V)$,
there exists an irreducible curve ${\cal C}$ of $V$ such that $z \in \overline{ f( {\cal C} ) }$
(see the corollary on page 262 of \cite{Kr}).
Therefore, we may assume that $V$ is an irreducible curve.
By normalizing $V$ and by Nagata's theorem (see \cite{Na}),
we may  suppose that $V$ is smooth
and that $W$ is complete.
Let $C$ be "the completion" of $V$, i.e.  a smooth
projective curve containing $V$ as an open subset.
Since $W$ is complete,
$f$ can be (uniquely) extended in a morphism
$f: \, C \to W$.
We have $\overline{ f(V) }= f(C)$,
so that it is enough to show that for any point $x \in C $,
there exists $\varphi \in V \! \! \left( \rule{0mm}{4mm}  \C ((t)) \right)$
such that $x = {\displaystyle \lim_{t \to 0} } \, \varphi (t)$.
We can  assume that $x \notin V$ because otherwise there
is nothing to do.
Finally, taking a well chosen affine neighborhood of $x$ in $C$, we
can suppose that $C$ is affine and that $V= C \setminus \{ x \}$.
Let ${\cal O}(C)$ be the algebra of regular functions on $C$,
let ${\cal O}_{C,x}$ be the local ring of $x$ on $C$
and let $\widehat{ {\cal O}_{C,x} }$ be its completion.
We have natural injections
${\cal O}(C) \hookrightarrow {\cal O}_{C,x} \hookrightarrow \widehat{ {\cal O}_{C,x} }$
and it is well-known that $\widehat{ {\cal O}_{C,x} } \simeq \C[[t]]$.
Let $\C (C) \hookrightarrow \C((t ))$ be the extension to  fields of fractions
of the map ${\cal O}(C) \hookrightarrow \C[[t]]$.
We have the commutative diagram:

\[ \xymatrix{
{\cal O}(C) \ar@{^{(}->}[d] \ar@{^{(}->}[r] & {\cal O}_{C,x} \ar@{^{(}->}[r] 
& \widehat{ {\cal O}_{C,x} } \ar[r]^{\thicksim} & \C [[ t]] \ar@{^{(}->}[dd] \\
{\cal O}(V)\ar@{^{(}->}[d] \ar@{^{(}->}[rrrd]^{\varphi ^*}         \\
\C (C)\ar@{^{(}->}[rrr] & & & \C ((t)) }
\]

\noindent where $\varphi ^*: \, {\cal O}(V) \to \C ((t))$ is the algebra morphism
corresponding to the point $\varphi : {\rm Spec}\, \C((t)) \to V$ which
we were looking for.

\vskip1mm

(ii) $\Longrightarrow$ (i). It is well-known. \cqfd

\vskip5mm

\noindent {\bf Remark.}  Note the analogy with the metric case
where $w \in \overline{ f(V) }$ if and only if there exists
a sequence $(v_n)_{n \geq 1}$ of $V$ such that
${\displaystyle w = \lim_{n \to + \infty} f(v_n) }$.

\vskip5mm

\noindent {\bf Proof of lemma 3.1.}
Assume that $\gamma := ( \alpha X, \beta Y) \in  \overline{ {\cal C}(f)_{\leq \, m} }$.

If $\Omega:= \{ a^k b^l, \, 0 \leq k+l \leq m \}$,
the proof of lemma 2.6 tells us that $\alpha, \beta  \in \Omega \subseteq  < a,b>$,
so that $   < \alpha , \beta > \,  \subseteq \,  <a,b>$.

\vskip3mm

Let us prove the reverse inclusion.
By proposition 2.3, ${\cal C}(f)_{\leq \, m}$ is included in the image
of the map 
${\cal G}_{\leq \, m} \to {\cal G}$,
$\varphi \mapsto  \varphi ^{-1} \circ f \circ \varphi $.
Using the above valuative criterion, we get
the existence of
$\varphi  \in {\cal G}_{\leq \, m} \left( \rule{0mm}{4mm} \C((t)) \right)$ such that
if $g:= \varphi ^{-1} \circ f \circ \varphi 
\in {\cal G} \left( \rule{0mm}{4mm} \C((t)) \right)$,
then
$\gamma = {\displaystyle \lim_{t \to 0} } \, g_t$.
We have $g_t^*= \varphi_t^* \circ f^* \circ  (\varphi_t^*)^{-1}$ as linear endomorphisms
of the $\C((t))$-vector space $\C((t)) [X,Y]$.
Therefore $u_{k,\, l}:=\varphi_t^*(X^kY^l)$ is an eigenvector of $g_t^*$
associated with the eigenvalue $a^kb^l$.
Let $m \in \Z$ be such that  $v_{k,\, l}:=t^mu_{k,\, l}$
admits a nonzero limit $\overline{ v_{k,\, l} }$ when $t$ goes to zero.
We have $g_t^*( v_{k,\, l})= a^kb^l v_{k,\, l}$ and setting $t=0$, we get
$\gamma^*( \overline{ v_{k,\, l} } ) =a^kb^l \overline{ v_{k,\, l} }$.
Hence $a^kb^l$
is an eigenvalue of $\gamma ^*$, so that $a^kb^l \in < \alpha , \beta >$. \cqfd

\vskip5mm

{\bf 2. Topological lemmas.}

\vskip5mm

\noindent {\bf Lemma 3.2.} Let $f=(aX,bY) \in {\cal G}$.
If $( \alpha X, \beta Y) \in   \overline{ {\cal C}(f)_{\leq \, m} }$
with $\alpha, \beta \neq 1$, then $\{ \alpha , \beta \} = \{ a,b \}$.

\vskip5mm

\noindent {\bf Proof.}

\noindent \underline{Claim.} For any $\varepsilon >0$
there exists a $C^0$-neighborhood $U$
of $\gamma := ( \alpha X, \beta Y)$ in ${\cal E}_{\leq \, m}$
such that any  $g \in U$ admits a fixed point
in $B_{\varepsilon}$.

\vskip1mm

Indeed, there exists an $\eta >0$  such that
$ B_{\eta} \subseteq (\gamma - {\rm id})(B_{\varepsilon}) $,
so that there exists a $C^0$-neighborhood $U$
of $\gamma $  such that  any  $g \in U$ satisfies $0 \in (g - {\rm id})(B_{\varepsilon})$.

\vskip2mm

Let $(g_n)_{n \geq 1}$ be  a sequence of ${\cal C}(f)_{\leq \, m}$
such that $\gamma = {\displaystyle \lim_{n \to \infty} }g_n$
for the $C^1$-topology.
By the  claim, there exists a sequence $(\xi_n)_{n \geq 1}$
of points of $\A^2$ such that $g_n(\xi_n)=\xi_n$
and ${\displaystyle \lim_{n \to \infty} }\xi_n=0$.
Therefore, we have $\gamma '(0) = {\displaystyle \lim_{n \to \infty} }g_n'(\xi_n)$
for the usual topology of $M_2(\C)$.
Since ${\rm Tr}\, \gamma '(0)= \alpha + \beta$
and ${\rm Tr}\,g_n'(\xi_n) = a+b$, we get $\alpha + \beta = a+b$.
But  $\alpha  \beta = ab$ (using the Jacobian), so that
$\{ \alpha , \beta \} = \{ a,b \}$. \cqfd

\vskip5mm

We will admit the following convexity lemma.

\vskip5mm

\noindent {\bf Lemma 3.3.} If $B'$ is a closed ball in an euclidian
space, there exists a $C^2$-neighborhood  of the
identity map on the space such that for any $g$ in this neighborhood, $g(B')$ is convex.

\vskip5mm

\noindent {\bf Remark.} Let 
$B':=\{ \rho e^{ i \theta }, \, \theta \in \R,  \, 0 \leq \rho \leq 1  \}$
be the unit disc in $\C$.
If $g$ is "near" the identity  for the $C^2$-topology, then
we will have
$g(B')= \{ \rho e^{ i \theta },  \, \theta \in \R,  \, 0 \leq \rho \leq r( \theta ) \}$
where $r : \R \to \R$ is a $2 \pi$ -periodic map which
is "near" the map $s \equiv 1$ for the $C^2$-topology.
The curvature of the parametrized curve $\theta \mapsto r( \theta) e^{i \theta }$
at the point $\theta $ is well-known to be
${\displaystyle C= \frac{r^2 + 2 \, {r'}\,^2-r \, r''}{(r^2+{r'}\,^2)^{\frac{3}{2}}} }$.
If $r$ is "near"  $s$ for the $C^2$-topology,
it is clear that $C>0$
at each point, showing that $g(B')$ is convex.

\vskip5mm

\noindent {\bf Lemma 3.4.} If $f$ is a finite-order automorphism,
${\cal C}(f)$ is closed in ${\cal G}$.

\vskip5mm

\noindent {\bf Proof.} We may assume that $f=(aX,bY)$
where $a^q=b^q=1$ for some $q \geq 1$.
It is enough to show that 
if $\gamma=( \alpha X, \beta Y) \in \overline{ {\cal C}(f)_{\leq \, m} }$
for some $m$,
then $\{ \alpha, \beta \}= \{ a,b \}$.

We begin to note that $g^q= {\rm id}$
for any  $g \in {\cal C}(f)$.

\vskip2mm

\noindent \underline{Claim.} For any $\varepsilon >0$
there exists a $C^2$-neighborhood $U$ of $\gamma $ in ${\cal E}_{\leq \, m}$
such that if $g \in U$ and  $g^q={\rm id}$,
then $g$ admits a fixed point in $B'_{\varepsilon}$.

Let us note that $\gamma ( B'_{\varepsilon}  ) = B'_{\varepsilon}$.
It is enough to take  for $U$ a $C^2$-neighborhood  of
$\gamma$ such that for any $g \in U$
and any $0 \leq k < q$, $g^k( B'_{\varepsilon} )$
is a convex set containing the origin.
Indeed, if $g \in U$ and  $g^q = {\rm id}$, then
$K:= {\displaystyle \bigcap_{0 \, \leq \,  k \,  < \, q} } g^k( B'_{\varepsilon } )$
is a non-empty compact convex set such that $g(K)=K$.
By Brouwer fixed point theorem, $g$ admits a fixed point in
$K \subseteq B'_{\varepsilon}$ and the claim is proved.

\vskip2mm

We finish the proof exactly as in lemma 3.2. \cqfd

\vskip10mm

{\bf 3. The proof.}

\vskip5mm

($\Longrightarrow$)
Thanks to proposition 2.3 it is enough to show that if
$f=(aX,bY) \in {\cal G}$, then
$ {\cal C}(f)$ is closed in ${\cal G}$.
Thanks to lemma 2.6 it is enough to show that if 
$\gamma= ( \alpha X, \beta Y) \in \overline{ {\cal C}(f)_{\leq \, m} }$
for some $m$, then $\{ \alpha, \beta \}=\{ a,b \}$.

\vskip4mm

\noindent \underline{First case.} $\alpha, \beta \neq 1$.

We conclude by lemma 3.2.

\vskip4mm

\noindent \underline{Second case.} $\alpha$ or  $\beta = 1$.
We can assume that $\alpha =1$.

Since ${\rm Jac}\, \gamma= {\rm Jac}\, f$, we have $\beta = ab$.
But $<a,b>= < \beta >$ by lemma 3.1, so that there exist
$k,l \geq 0$  such that $a = \beta ^k$, $b= \beta ^l$.

\vskip2mm

\underline{First  subcase.} $\beta$ is not a root of unity.

The equality $\beta =ab$ gives us
$\beta = \beta ^{k+l}$, so that $1 = k + l$.
We get $\{k,l\}=\{ 0,1 \}$, so that $\{a,b \}= \{ 1, \beta \}= \{ \alpha ,\beta \}$.

\vskip2mm

\underline{Second  subcase.} $\beta$ is  a root of unity.

It is clear that $a,b$ are also roots of unity.
Therefore, $f$ is a finite-order automorphism
and we conclude by lemma 3.4.

\vskip5mm

($\Longleftarrow$) Let $f$ be any polynomial automorphism.
We want to show that  $\overline{ {\cal C}(f) }$
contains a semisimple polynomial automorphism.
It is  sufficient to show that it contains  a linear automorphism.
Indeed,  in the linear group
it is well-known that any conjugacy class contains in its closure
a (linear) semisimple automorphism.

\vskip2mm

\noindent \underline{First case.} $f$ is triangularizable.

We can assume that $f=(aX+p(Y),bY+c)$.
If $l_t:=(tX,Y)$ and $r_t:=(X,tY) \in {\cal G}$ for  $t \in \C^*$,
we have ${\displaystyle \lim_{t \to 0} }\, l_t \circ f \circ (l_t)^{-1} \, = \, (aX,bY+c)$.
Therefore, $u:= (aX,bY+c) \in \overline{ {\cal C}(f) }$.
But $r_t \circ u \circ (r_t)^{-1} \in \overline{ {\cal C}(f) }$ for any $t \neq 0$
and ${\displaystyle \lim_{t \to 0} }\, r_t \circ u \circ (r_t)^{-1}=(aX,bY)$.

\vskip2mm

\noindent \underline{Second case.} $f$ is not triangularizable.

We can assume that $f$ is cyclically reduced 
of degree $d \geq 2$.
By theorem 3.1 of \cite{F-M}, $f$ has exactly $d$ fixed points (counting the multiplicities).
In particular, it has a fixed point and by conjugating we can assume that
it fixes the origin.
Therefore, if  $h_t:= (tX,tY) \in {\cal G}$ for $t \neq 0$,
then ${\displaystyle \lim_{t \to 0}\, (h_t)^{-1} \circ f \circ h_t }$ is equal to the linear part
of $f$.  \cqfd

\vskip5mm 

\noindent {\bf Acknowledgements.} We would like to thank
Hanspeter Kraft for  suggesting to us this question
and  the staff of the Mathematisches Forschungsinstitut
Oberwolfach for its hospitality.


\vskip5mm 






\footnotesize
 
\begin{thebibliography}{99}



\bibitem{B-C-W}  H. Bass, E. Connell, D. Wright, The Jacobian Conjecture:
Reduction of Degree and Formal Expansion of the Inverse,
Bull. of the A.M.S., 7 (1982), 287-330.



\bibitem{CMP} L. Cerlienco, M. Mignotte, F. Piras, Suites r\'ecurrentes lin\'eaires,
propri\'et\'es alg\'ebriques et arithm\'etiques,
L'Enseignement Math\'ematique 33 (1987), 67-108.




\bibitem{E} D. Eisenbud, Commutative algebra. With a view toward algebraic geometry.
Graduate Texts in Mathematics, 150, Springer-Verlag, New York, 1995.



\bibitem{F-M} S. Friedland, J. Milnor, Dynamical properties
of plane polynomial automorphisms, Ergod. Th $\&$ Dyn. Syst. 9 (1989),
67-99.


\bibitem{Fu1} J.-P. Furter, On the degree of iterates of automorphisms of the affine plane,
Manuscripta Mathematica, 98 (1999), 183-193.


\bibitem{Fu2} J.-P. Furter, Quasi-Locally Finite Polynomial Endomorphisms,
submitted for publication, available at
http://www.univ-lr.fr/labo/lmca/publications/06-07/06-07.pdf.




\bibitem{Fu-Ma} J.-P. Furter, S. Maubach,
Locally Finite Polynomial Endomorphisms,
Journal of Pure and Applied Algebra 211 (2007), no. 2, 445-458.



\bibitem{Gr}  A. Grothendieck, J. Dieudonn\'e,
El\'ements de g\'eom\'etrie alg\'ebrique, II Etude globale
\'el\'ementaire de quelques classes de morphismes, Publ. Math.
IHES 8 (1961).




\bibitem{Ha} R. Hartshorne,  Algebraic geometry,
Graduate Texts in Mathematics, No. 52,
Springer-Verlag, New York-Heidelberg, 1977.



\bibitem{He} G. Hermann, 
Die Frage der endlich vielen Schritte in der Theorie der Polynomideale,
Math. Ann. 95 (1926), no. 1, 736-788.


\bibitem{Ju} H. W. E. Jung, \"{U}ber ganze birationale
Transformationen der Ebene, J. Reine Angew. Math. 184 (1942),
161-174.

\bibitem{Kr} H. Kraft,
Geometrische Methoden in der Invariantentheorie. (German) [Geometrical methods in invariant theory] 
Aspects of Mathematics, D1. Friedr. Vieweg \& Sohn, Braunschweig, 1984.



\bibitem{Ku} W. van der Kulk, On polynomial rings in two
variables, Nieuw. Arch. Wisk. (3) 1 (1953), 33-41.


\bibitem{M-M} 
E. W. Mayr, A.R. Meyer,
The complexity of the word problems for commutative semigroups and polynomial ideals,
Adv. in Math. 46 (1982), no. 3, 305-329.


\bibitem{Mu} D. Mumford,
The red book of varieties and schemes.
Second, expanded edition.
Lecture Notes in Mathematics,
Springer-Verlag, Berlin, 1999.



\bibitem{M-F-K} D. Mumford, J. Fogarty, F. Kirwan,
Geometric invariant theory. Third edition.
Ergebnisse der Mathematik und ihrer Grenzgebiete (2)
[Results in Mathematics and Related Areas (2)], 34.
Springer-Verlag, Berlin, 1994.



\bibitem{Na} M. Nagata, Imbedding of an abstract variety
in a complete variety, J. Math. Kyoto Univ. 2 (1962), 1-10.



\bibitem{Sei} A. Seidenberg, 
Constructions in algebra,
Trans. Amer. Math. Soc. 197 (1974), 273-313.


\bibitem{Se} J.-P. Serre, Trees, Springer Verlag, Berlin, Heidelberg, New York, 1980.




\bibitem{Sha1} I. R. Shafarevich, On some infinite-dimensional groups,
Rend. Mat. e Appl. (5) 25 (1966), no. 1-2, 208-212.



\bibitem{Sha2} I. R. Shafarevich, On some infinite-dimensional
groups II, Math. USSR Izv., 18 (1982), 214-226.




\end{thebibliography}
\end{document}